\documentclass[12pt]{article}
\def\supt{\sup_{0\leq t\leq 1}}

\newcommand{\has}{\hbox{ a.s.}}
\newcommand{\qas}{\hbox{ ~a.s.}}
\def\dto{\displaystyle\mathop{\rightarrow_d}}
\def\diag{\hbox{ diag}}

\def\noo{n\to\infty}

\def\deq{\stackrel{{d}}{\,=\,}}
\newcommand{\IR}{\mathbb R}
\newcommand{\N}{\mathbf N}
\newcommand{\bZ}{\mathbf Z}
\newcommand{\bW}{\mathbf W}
\newcommand{\IZ}{\mathbb Z}
\newcommand{\IW}{\mathbb W}

\newcommand{\Cs}{Cs\"org\H{o}}
\newcommand{\W}{\mathbf W}
\newcommand{\Z}{\mathbf Z}
\newcommand{\cD}{\mathcal D}
\newcommand{\cE}{\mathcal E}
\newcommand{\cS}{\mathcal S}
\newcommand{\ep}{\varepsilon}
\newcommand{\beq}{\begin{equation}}
\newcommand{\eeq}{\end{equation}}

\newtheorem{thm}{Theorem}

\usepackage[english]{babel}
        \usepackage{latexsym}
\usepackage{amssymb}

\begin{document}

\vskip.5cm\noindent
{\Large \bf Two-dimensional Anisotropic Random Walks:} 

\noindent
{\Large \bf Fixed {\sl versus}
Random Column Configurations for} 

\noindent
{\Large \bf Transport Phenomena}

\vskip.9cm\noindent
{\bf Endre Cs\'aki}

\noindent
A. R\'enyi Institute of Mathematics, Hung. Acad. Sci., Budapest, P.O.B.
127, H-1364, Hungary. E-mail address: csaki.endre@renyi.mta.hu

\vskip.3cm\noindent
{\bf Mikl\'os Cs\"org\H o}

\noindent
School of Mathematics and Statistics, Carleton University, 1125 Colonel By
Drive, Ottawa, Ontario, Canada K1S 5B6. E-mail address: 
mcsorgo@math.carleton.ca

\vskip.3cm\noindent
{\bf Ant\'onia F\"oldes}

\noindent
Department of Mathematics, College of Staten Islands, CUNY, 2800 Victory 
Blvd., Staten Island, New York 10314, U.S.A. E-mail address: 
Antonia.Foldes@csi.cuny.edu

\vskip.3cm\noindent
{\bf P\'al R\'ev\'esz}

\noindent
Institut f\"ur Statistik und Wahrscheinlichkeitstheorie, Technische 
Universit\"at Wien, Wiedner Hauptstrasse 8-10/107 A-1040 Vienna, Austria.

\noindent
E-mail address: revesz.paul@renyi.mta.hu

\vskip.5cm\noindent
{\bf Abstract}
We consider random walks on the square lattice 
of the plane along the lines of Heyde (1982, 1993) and den Hollander (1994), 
whose studies have in part been inspired by the so-called transport phenomena 
of statistical physics. Two-dimensional anisotropic random walks with  
anisotropic density conditions {\it \'a \,la } Heyde (1982, 1983) yield 
{\it fixed column configurations} and nearest-neighbour random walks in a 
random environment on the square lattice of the plane as in den Hollander 
(1994) result in {\it  random column configurations}. In both cases  we 
conclude simultaneous weak Donsker and strong Strassen type invariance 
principles in terms of appropriately constructed anisotropic Brownian motions 
on the plane, with self-contained proofs in both cases. The style of 
presentation throughout will be that of a semi-expository survey of related 
results in a historical context.


\noindent
anisotropic Brownian motions; in probability and almost sure stochastic
approximations

\section{Introduction}\label{sect1}

Motivated by, and along the lines of, Heyde (1982, 1993) and den Hollander 
(1994), in this paper we consider random walks on the square lattice 
of the plane whose studies have in part been inspired by the so-called 
transport phenomena of statistical physics (cf., e.g., {\bf 1 Introduction} 
in Heyde (1993) and {\bf 1.4 History} of den Hollander (1994), and their 
references). 

Two-dimensional anisotropic random walks with asymptotic density conditions
{\it \'a \,la } Heyde (1982, 1983) yield {\it fixed column configurations }
(cf. Section 1.1) and nearest-neighbour random walks in a random environment 
on the square lattice of the plane as in den Hollander (1994) result in 
{\it random column configurations} (cf. Section 1.2). In both cases, in the 
same respective sections, we conclude simultaneous weak Donsker and strong 
Strassen type invariance principles in terms of appropriately constructed 
anisotropic Brownian  motions on the plane via making use of the initial  
approach of the above papers in hand. 
	
Thus, in Section 1.1, via the initial approach of Heyde (1992, 1993), 
we arrive at the respective weak and strong fixed column embedding conclusions 
as in Theorem 1 and Theorem 2, holding simultaneously on the same  probability 
space, while in Section 1.2 we conclude the respective corresponding  
simultaneous weak and strong random column conclusions of Theorem C and 
Theorem 3, via the initial approach of den Hollander (1994).

In view of this syncronized parallel  relationship we may say that, to the 
extent of simultaneously having a weak Donsker and strong Strassen type 
asymptotic  behaviour, under their  respective initial  conditions the two 
random walks in hand behave similarly. More precisely, in addition to 
extending the results of Heyde (1982, 1993) on the weak Donsker and strong 
Strassen type asymptotic path behaviour of a 2-dimensional anisotropic random 
walk with transition probabilities as in (1.1) under his asymptotic density 
conditions as in (1.2), the latter with $\eta=0$ in the case of Theorem 1,  
in this exposition we have also succeeded in showing that, to the extent of 
weak Donsker and strong Strassen type asymptotic path behaviour, a 
nearest-neighbour 2-dimensional random walk with probability law $\mu$ as in 
the ($\star$) den Hollander (1994) condition of (1.17) and transition 
probabilities as in (1.19)  behaves like a 2-dimensional anisotropic random 
walk with anisotropy determined only by $q$ as in (1.17) and not by any other 
parameters of $\mu$. This in turn  also  amounts to the main  observation of 
this exposition (cf. Observation in Section 1.2).
	
The corresponding respective self contained proofs are given in Sections 2.1 
and 2.2. The style of presentation throughout will be that of a 
semi-expository survey of related  results in a historical context.
	
We also wish to note that, using an alternative construction of the random 
walk $\{X_n,Y_n\},\,\,n=1,2,....,$ on ${\mathbb Z^2}$ with horizontal step 
transition probabilities $1/2-p_j$ and vertical step transition probabilities 
$p_j,$ i.e., an inverted "upside down"  version of transion probabilities as 
in (1) of Heyde (1993) (cf. (1.1) in this exposition), and assuming the 
asymptotic density condition of Heyde (1982) for the transition probabilities 
as spelled out in (1.2), Cs\'aki {\it et al}. (2013) established a joint 
strong approximation result for $\{X_n,Y_n\}$ by two independent standard 
Wiener processes as in their Theorem 1.1. Mutatis mutandis, an alternative 
equivalent version of our Theorem 2 per se could also be proved via making 
use of the latter strong approximation result. On the other hand, our Donsker  
type in probability approximation as in Theorem 1 is new, and so is also its 
relationship to our present way of establishing Theorem 2 (cf. Remark 2).
	
Next we note that our strong Strassen type Theorem 3, an extended version of  
(2.1) of Theorem 2 in den Hollander (1994), is new and so is also its 
relationship under the same conditions to Therem C, a reformulated version of 
Theorem 1 of den Hollander (1994) {\it \'a la} our Theorem 1 (cf. Remark 5).
	
In the context of nearest-neighbour random walks in random environment on the 
square lattice on the plane, den Hollander (1994) also deals  with other 
functions of the walk, like return probability and range. For anisotropic 
random  walks on the plane  Cs\'aki {\it et al}. (2013) also study their 
local times, recurrence and range under various conditions. We have not 
succeeded in furthering and synchronizing the result yielded by the two 
approaches in hand along the latter lines.

\subsection {Two-dimensional anisotropic random walks; Fixed column 
configurations with asymptotic density conditions}\label{ssect1.1}  

\noindent 
For the sake of introducing and motivating the study of two-dimenssional 
anisotropic random walks on the square lattice ${\mathbb Z^2}$ we quote from 
1.~Introduction of Heyde (1993):

\begin{quote}

"One area of application involves conductivity of various organic salts, such 
as tetrathiofuvalene (TTF)-tetracyanoquinodimethane (TCNQ), which show signs 
of superconductivity. They conduct strongly in one direction but not others. 
Indeed, the conductivity parallel to the structural axis is 100 times or more 
that perpendicular to it. The molecular forms of TCNQ and TTF are planar so 
that they can be easily stacked and, in fact, the structure of the TTF-TCNQ 
crystal is generally thought to consist of parallel columns of separately 
stacked TTF and TCNQ molecules.

Another related area of application concerns transport in physical systems 
which lack complete connectivity such as with fluid flow through packed 
columns, for example in gas absorption or distillation processes. In such 
situations gases or liquids predominantly flow through vertical cylinders 
filled with a random packing of inert objects such as ceramic rings.

The predominance of one dimension for the transport, say along rows, plus a 
possible lack of complete connectivity, can be modeled rather more generally 
by an anisotropic 2-dimensional random walk in which the transition mechanism 
depends only on the index of the column which is at present occupied. Thus, 
we consider a random walk which, if situated at a site on column $j$, moves 
with probability $p_j$ to either horizontal neighbour and with probability 
$\frac{1}{2} - p_j$ to either vertical neighbour at the next step."
\end{quote}

More formally, let $X_n$ and $Y_n$ denote the horizontal and vertical 
positions of the walk after $n$ steps, starting from $X_0=Y_0=0$, with 
transition probabilities as in (1) of Heyde (1993):
\begin{eqnarray}
P\{(X_{n+1},Y_{n+1}) &=& (j-1,k)|(X_n, Y_n)=(j,k)\}=p_j\nonumber\\
P\{(X_{n+1},Y_{n+1}) &=& (j+1,k)|(X_n, Y_n)=(j,k)\}=p_j\nonumber\\
P\{(X_{n+1},Y_{n+1}) &=& (j,k\! -\!1)|(X_n,Y_n)=(j,k)\}
= 1/2\! -\! p_j\nonumber\\
P\{(X_{n+1},Y_{n+1}) &=& (j,k +1)|(X_n,Y_n)=(j,k)\}= 1/2\! -\! p_j
\label{eq1.1}
\end{eqnarray}
for $(j,k) \in \IZ^2$, and $n=0,1,2,\ldots$, a so-called anisotropic 
2-dimensional random walk with possibly unequal {\bf symmetric} horizontal 
and vertical step transition probabilities that depend only on the index of 
the column which is at present occupied.

We assume throughout that $0\!<\!p_j \!\leq\! 1/2$ and 
min$_{j\in\IZ}\, p_j\! <\! 1/2$.

Thus, as in Heyde (1993), we consider a random walk which, if situated at
a site on column $j$, moves with probability $p_j$ to either horizontal 
neighbour, and with probability $1/2-p_j$ to either vertical neighbour at
the next step. 

The case $p_j = 1/4$, $j=0, \pm 1, \pm 2, \ldots$ corresponds to a simple 
symmetric random walk on the plane. For studies of this case we refer to 
Erd\H{o}s and Taylor (1960), Dvoretzky and Erd\H{o}s (1951), R\'ev\'esz 
(2013) and the references in these works.

When $p_j = 1/2$ for some $j$, then the vertical line $x=j$ is missing, not 
connected, i.e., we have a so-called non-connective column.

If all $p_j = 1/2$, then the random walk takes place on the $x$ axis, i.e., 
$(X_n,Y_n)$ reduces to a simple symmetric random walk on the real line. This 
case is excluded from this investigation. Hence the assumption that 
min$_{j\in\IZ} p_j < 1/2$. 

When $p_j = 1/2$, $j=\pm 1, \pm 2, \ldots$, but $p_0=1/4$, then all the 
vertical lines $x=j=\pm 1, \pm 2, \ldots$ are missing, except that of $x=0$, 
i.e., the $y$ axis. Thus we get what could be called a random walk on a 
``horizontal comb'', whose vertical version was studied by Weiss and Havlin 
(1986), Bertacchi and Zucca (2003), Bertacchi (2006), Cs\'aki {\it et al.} 
(2009, 2011).

Back to the random walk (1.1), even though the transition probabilities depend 
on the position of the first coordinate of ${\mathbf Z}_n := (X_n,Y_n)$, 
$n\in {\mathbf N}$, both $X_n$ and $Y_n$ will be seen to behave like a simple 
symmetric random walk on $\IZ$, except for a random time delay, and 
independently of each other.

Heyde (1982) assumed the following asymptotic density condition for the 
transition probabilities of (1.1):

\beq
k^{-1} \sum^k_{j=1} p_j^{-1}=2\gamma + o(k^{-\eta}),\quad \ \ \ \ \ 
~~ k^{-1} \sum^k_{j=1} p_{-j}^{-1} = 2\gamma + o(k^{-\eta})
\label{eq1.2}
\eeq
as $k\to\infty$, for some constants $1<\gamma<\infty$ and $1/2 < \eta <\infty$.

Under condition (1.2), for the first coordinate of $\bZ_n$, Heyde (1982) 
concluded his Theorem~1 that reads as follows. 

\noindent{\bf Theorem A}~
{\sl On an appropriate probability space for $X_n$ there is a standard Wiener 
process $\{W(t),t\geq 0\}$ so that}
\beq
\gamma^{1/2}X_n = W(n(1+\varepsilon_n)) + 
O(n^{1/4}(\log n)^{1/2}(\log\log n)^{1/2})\ \has
\label{eq1.3}
\eeq
as $n\to\infty$, {\sl with} $\varepsilon_n\to 0\,$\ \ a.s. 

As to this theorem, Heyde (1982) writes:

\begin{quote}

``If strong additional assumptions are made about the column types, for 
example if there is a finite number of different types which occur with 
fixed periodicity, then the random variable $\ep_n$ can safely be removed 
and relegated to the error in (1.3). We shall not pursue this question here 
because the presence of $\epsilon_n$ does not cause undue complication in 
the extraction of specific results from (1.3) as evidenced by Corollary 1 
below.''

\end{quote}

The just mentioned Corollary 1 in Heyde (1982) reads as follows.

\noindent{\bf Corollary A} ~{\sl Under condition $(1.2)$, we have}
\begin{itemize}
\item[(i)] $\gamma^{1/2}n^{-1/2}X_n \dto N(0,1)$, ~~{as~ } $n\to \infty$,
\item[{(ii)}] $\limsup_{n\to\infty} (2n\gamma^{-1}\log\log n)^{-1/2} X_n = 
1 ~\hbox{ a.s. }$,
\item[{(iii)}] $\liminf_{n\to\infty} (2n\gamma^{-1}\log\log n)^{-1/2} X_n 
= - 1 ~\hbox{ a.s. }$
\end{itemize}

One of the aims of this exposition is to show that no additional conditions 
are needed for removing and relegating the random variable $\ep_n$ to the 
error in (1.3) (cf. (1.13) of our Theorem 2, Corollary 2 and our conclusion 
right after the latter corollary). Let 
$\sigma_0 = 0 < \sigma_1 < \sigma_2 < \cdots$ be the successive times at 
which the values of  $X_i-X_{i-1}$, $i=1,2,\ldots$, are nonzero, and put 
$S_1(k)=X_{\sigma_k}$. By the assumed symmetry of the transition probabilities 
in (1.1), $S_1(k)$ is a simple symmetric random walk on $\IZ$.

In view of this and (1.2), Heyde (1982) concludes his Lemma 1: 
$\sum_1^\infty n^{-2} Ep^{-2}_{S_n} < \infty$ {\sl and hence} 
$\sum_1^\infty n^{-2}p^{-2}_{S_n} < \infty$ \ a.s.  This, in turn, leads 
to the principal result that is needed for establishing his Theorem 1, i.e.,
Theorem A as above, namely to Theorem 2 of Heyde (1982):

\beq
n^{-1}\sigma_n \to \gamma \hbox{ a.s., ~as} ~n\to\infty.
\label{eq1.4}
\eeq

Also, $X_n = X_{\sigma_k}$ for $\sigma_k \leq n < \sigma_{k+1}$.  As in the 
proof of Theorem 1 of Heyde (1982), for $n$ fixed let
$$
\sigma_{k(n)} := \max\left[j:1 \leq j\leq n, \, X_j \neq X_{j-1}\right].
$$
Then
\beq
X_n = X_{\sigma_{k(n)}} = S_1(k(n))
\label{eq1.5}
\eeq
is the horizontal position of the walk ${\mathbf Z}_n = (X_n,Y_n)$ after 
$k(n)$ horizontal steps in the first $n$ steps of ${\mathbf Z}_n$ and, using 
(1.4), under condition (1.2) we have
\beq
n^{-1}k(n) \to \gamma^{-1} \hbox{ a.s., ~~as } n\to\infty
\label{eq1.6}
\eeq
(cf. Proof of Theorem 1 of Heyde (1982)). The latter conclusion in combination 
with Strassen's invariance principle (cf. Strassen (1967)) as in Heyde (1982), 
results in his Theorem 1 (cf. Theorem A above).

Clearly, in view of (1.5), $\ell(n) := n-k(n)$ is the number of vertical steps 
in the first $n$ steps of $\bZ_n$ and, as a consequence of (1.6), under 
condition (1.2) we also have  
\beq
n^{-1}\ell(n) \to 1-\gamma^{-1} \hbox{ ~a.s., ~as } ~n\to\infty.
\label{eq1.7}
\eeq

Under the weaker condition that $\eta = 0$ in (1.2), Heyde (1993) concluded 
(1.6) and (1.7) in probability, i.e., that as $n\to\infty$, we have
\beq
n^{-1}k(n) \to \gamma^{-1} \hbox{ ~and~ } n^{-1}\ell(n) 
\to 1-\gamma^{-1} \hbox{ ~in probability}.
\label{eq1.8}
\eeq

Using this together with the Cram\'er-Wold device in combination with the 
martingale limit result of Theorem 3.2 in Hall and Heyde (1980), Heyde (1993) 
established the following asymptotic joint distribution for the random walk 
${\mathbf Z}_n = (X_n,Y_n)$.

\noindent{\bf Theorem B} ~{\sl Assume $(1.2)$ with $\eta = 0$. Then}, 
as $n\to\infty$,
$$
n^{-1/2}{\mathbf Z}_n = 
\left(\frac{X_n}{n^{1/2}},\, \frac{Y_n}{n^{1/2}}\right) 
\dto \left(\gamma^{-1/2}N_1,(1-\gamma^{-1})^{1/2}N_2\right),
\leqno(\hbox{a})
$$
{\sl where $N_1$ and $N_2$ are independent standard normal random variables, 
and}
$$
EX_n^2 \sim \gamma^{-1}n \hbox{ ~and~ } EY^2_n \sim (1-\gamma^{-1})n.
\leqno(\hbox{b})
$$
The symbol $\sim$ means that the ratio of the two sides tends to 1 as 
$n\to\infty$. Here, and throughout, $\dto$ stands for convergence in 
distribution.

Thus, as noted by Heyde (1993), the asymptotic behaviour of 
$n^{-1/2}(X_n,Y_n)$ depends only on the macroscopic properties of the medium 
through which the walk takes place, i.e., of the $\{p_j, \, j\in \IZ\}$ as in 
(1.2), with $\eta = 0$ in this particular case. In the same paper Heyde also 
notes that this is of considerable practical significance since, although many 
materials are heterogeneous on a microscopic scale, they are essentially 
homogeneous on a macroscopic or laboratory scale. In conclusion Heyde (1993) 
writes  (cf.\ {\bf  3.~Final Remark}): 

\begin{quote}

``The so called {\bf dimensional anisotropy} is $(1-\gamma^{-1})/\gamma^{-1}$ 
and for a material like TTF-TCNQ we expect this to be, say $10^{-2}$.  It 
should be noted that this would be achieved, for example, if only every 100th 
column were connective. We would have $p_j = \frac{1}{2}$ for $j$ a 
non-connective column and $p_j = \frac{1}{4}$ for $j$ a connective column."

\end{quote}

\noindent{\bf Remark 1} \label{rem1}
~We note that condition (1.2) clearly holds true if $\{p_j\}_{j\in\IZ}$ is a 
periodic sequence. Namely, in this case, for a given positive integer 
$L\ge 1$, $p_{j+L}=p_j$ for all $j\in\IZ$, and for $i=0,1,\ldots,$ we have
\beq
\frac{1}{L} \sum^{L-1}_{j=0} \frac{1}{p_{j+iL}} = 
\frac{1}{L} \sum^{L-1}_{j=0} \frac{1}{p_j} = 2\gamma,
\label{eq1.9}
\eeq
with some constant $1<\gamma<\infty$, as in (1.2). A particular periodic case, 
the so-called uniform case, when $p_j = 1/4$ if $|j|\equiv 0$ (mod L) and 
$p_j=1/2$ otherwise yields (1.9) with $\gamma = (L+1)/L$. This uniform 
periodic case may serve as a model for describing dimensional anisotropies as 
in the just quoted example of Heyde (1993) via 
$(1-\gamma^{-1})/\gamma^{-1} = \gamma -1 = 1/L$.  

Cs\'aki {\sl et al.} (2013) study the asymptotic behaviour of the anisotropic 
random walk $\{X_n,Y_n\}$ under both conditions (1.2) and (1.9). 

In order to formulate our results here, consider $D([0,\infty),\IR^2)$, the 
space of $\IR^2$ valued {\sl c\`adl\`ag\/} functions on $[0,\infty)$. For 
functions $(f(t),g(t)) = \Big((f_1(t),f_2(t))$, $(g_1(t),g_2(t))\Big)$ in this 
function space, define for all fixed $T>0$
\beq
\Delta = \Delta_T (f,g) := 
\sup_{0\leq t\leq T} \Vert (f_1(t)-g_1(t)),(f_2(t)-g_2(t))\Vert,
\label{eq1.10}
\eeq
where $\Vert \cdot \Vert$ is a norm in $\IR^2$, usually the 
$\Vert \cdot \Vert_p$ norm with $p=1$ or 2 in our case.

Based on Heyde's result as in (1.8) that is implied by condition (1.2) with 
$\eta = 0$, our first new conclusion is a Donsker type weak invariance 
principle, an extension of Heyde's bivariate central limit theorem (CLT) as 
in (a) of Theorem B, that reads as follows.

\begin{thm}\label{thm1} 
With $X_0 = Y_0 = 0$, let $\{X_n,Y_n,n\geq 0\}$ be the horizontal and vertical 
positions after $n$ steps of the random walk on $\IZ^2$ with transition 
probabilities as in $(1.1)$, and assume that condition $(1.2)$ holds true 
with $\eta = 0$.  Then, on an appropriate probability space for this random 
walk ${\mathbf Z}_n := \{X_n,Y_n\}$ on $\IZ^2$, one can construct two 
independent standard Wiener processes $\{W_1(t),t\geq 0\}$, 
$\{W_2(t),t\geq 0\}$  so that, as $n\to\infty$, with
\beq
\Big\{ {\mathbf W}_n(t), t\geq 0\Big\}_{n\geq 0} := \left\{
\frac{W_1(nt\gamma^{-1})}{n^{1/2}}, 
\frac{W_2(nt(1-\gamma^{-1}))}{n^{1/2}}, t\geq 0 \right\}_{n\geq 0}
\label{eq1.11}
\eeq
we have for all fixed $T>0$
\begin{eqnarray}
\lefteqn{\sup_{0\leq t\leq T} \Big\Vert n^{-1/2} {\mathbf Z}_{[nt]} - 
{\mathbf W}_n(t)\Big\Vert}\nonumber\\[2ex]
&=& \sup_{0\leq t\leq T} \left\Vert
\frac{X_{[nt]} - W_1(nt\gamma^{-1})}{n^{1/2}},
\frac{Y_{[nt]} - W_2(nt(1-\gamma^{-1}))}{n^{1/2}} \right\Vert \nonumber\\[2ex]
&=& o_P(1).
\label{eq1.12}
\end{eqnarray}
\end{thm}

Moreover, based on (1.6) instead of (1.8) that is implied by assuming (1.2) 
as postulated with $1/2 < \eta < \infty$, we arrive at a strong Strassen type 
invariance principle for ${\mathbf Z}_n := \{X_n,Y_n\}$ on $\IZ^2$ as follows.

\begin{thm}\label{thm2} 
With $X_0=Y_0=0$, let $\{ X_n,Y_n, \, n\geq 0\}$ be the horizontal and 
vertical positions after $n$ steps of the random walk on the integer lattice 
$\IZ^2$ with transition probabilities as in $(1.1)$ and assume that  
condition $(1.2)$ holds true. Then, on an appropriate probability space for 
this random walk ${\mathbf Z}_n := \{X_n,Y_n\}$ on $\IZ^2$, one can construct 
two independent standard Wiener processes $\{ W_1(t),t\!\geq\! 0\},$ 
$\{ W_2(t),\, t\geq 0\}$ so that, as $n\to\infty$, with 
$\left\{ {\W}_n(t), \,  t\geq 0 \right\}_{n\geq 0}$ as in $(1.11)$, we have 
\begin{eqnarray}
\lefteqn{\sup_{0\leq t\leq 1} \Big\Vert(2n\log\log  n)^{-1/2} 
{\mathbf Z}_{[nt]} -(2\log\log n)^{-1/2} {\mathbf W}_n(t)\Big\Vert}
\nonumber\\[2ex]
&=& \sup_{0\leq t\leq 1} \left\Vert \left(
\frac{X_{[nt]} - W_1(nt\gamma^{-1})}{(2n\log\log n)^{1/2}},
\frac{Y_{[nt]} - W_2(nt(1-\gamma^{-1}))}{(2n\log\log n)^{1/2}} \right) 
\right\Vert \nonumber\\[2ex]
&=& o(1) \ \rm{ a.s.}
\label{eq1.13}
\end{eqnarray}
\end{thm}

\noindent{\bf Remark 2}
\label{rem2}
{~It will be seen when proving Theorems 1 and 2 that the respective 
conclusions of (1.12) and (1.13) hold simultaneously on the same probability 
space in terms of the there constructed two independent standard Wiener 
processes $W_1$ and $W_2$. Thus $\W_n$ in (1.13) coincides with that of (1.11) 
not only by notation, but also by construction.}

Define now the measurable space $(D([0,\infty),\IR^2),\cD)$, where $\cD$ is 
the $\sigma$-field generated by the collection of all $\Delta$-open balls 
for all $T>0$ of the function space $D([0,\infty),\IR^2)$, where a ball is
a subset of $D([0,\infty),\IR^2)$ of the form
$$
\{(f_1,f_2):\, \Delta_T((f_1,f_2),(g_1,g_2))<r\}
$$
for some $(g_1,g_2)\in D([0,\infty),\IR^2)$ and $r>0$, with $\Delta_T$ as
in (1.10).

As a consequence of (1.12) of Theorem 1, on taking $T=1$, we may, for example, 
conclude weak convergence for the {\sl c\`adl\`ag\/} process 
$n^{-1/2}$ $\Z_{[n\cdot]}$ in  $D([0,1],\IR^2)$ in terms of the following 
functional convergence in distribution statement.

\noindent{\bf Corollary 1} ~{\sl Under the conditions of\/} 
{\rm Theorem 1}, {\sl as $n\to\infty$, we have}
\begin{eqnarray}
h(n^{-1/2}\Z_{[nt]}) &=& h\left(\frac{X_{[nt]}}{n^{1/2}}, 
\frac{Y_{[nt]}}{n^{1/2}}\right) \dto 
h\Big({\mathbf W}(t)
\hbox{\rm diag}(\gamma^{-1/2},(1-\gamma^{-1})^{1/2})\Big)\nonumber\\
&\deq & h(W_1(t\gamma^{-1}), W_2(t(1-\gamma^{-1}))
\label{eq1.14}
\end{eqnarray}
{\sl for all $h:D([0,1], \IR^2) \to \IR^2$ that are $(D([0,1],\IR^2),\cD)$ 
measurable and $\Delta$-continuous, or $\Delta$-continuous except at points 
forming a set of measure zero on $(D([0,1],\IR^2),\cD)$ with respect to the 
Wiener measure $\IW$ of $({\mathbf W}(t),0\leq t\leq 1) :=( (W_1(t), W_2(t)), 
0\leq t\leq 1)$, where $W_1$ and $W_2$ are two independent standard Wiener 
processes, and $\deq$ stands for equality in distribution.}

The Brownian motion ${\mathbf W}(t)
\hbox{\rm diag}(\gamma^{-1/2},(1-\gamma^{-1})^{1/2})$ on $\IR^2$ with the 
indicated diffusion matrix ${\rm diag}(\cdot,\cdot)$ is an example of an 
anisotropic Brownian motion.

\noindent{\bf Remark 3}\label{rem3}
~In the context of Corollary 1, $\cD$ of the measurable space 
$(D([0,1],\IR^2),\cD)$ coincides with the $\sigma$-field generated by the 
$\Delta$-{\sl open balls\/} of the function space $D([0,1],\IR^2)$ equipped 
with the norm $\Delta$ as in (1.8) with $T=1$. Clearly, the Brownian motion 
(Wiener process) $(\W(t),0\leq t\leq 1):=(W_1(t),W_2(t),0\leq t\leq 1)$, where 
$W_1$ and $W_2$ are two independent standard Wiener processes, is a random 
element of  $(C([0,1],\IR^2),\cE)$, the measurable space of $\IR^2$-valued 
continuous functions on $[0,1]$, where $\cE$ is the Borel $\sigma$-field 
generated by the $\Delta$-{\sl open subsets} of the space of $\IR^2$-valued 
continuous functions $C([0,1],\IR^2)$. Since
\beq
C([0,1],\IR^2) \in D([0,1],\IR^2) \hbox{ ~and~ } \cE = C([0,1],\IR^2) \cap \cD,
\label{eq1.15}
\eeq
$\W(t)$ is also a random element of $(D([0,1],\IR^2),\cD)$, and we can extend 
the Wiener measure $\IW$ of $\{\W(t),0\leq t\leq 1\}$ on 
$(C([0,1],\IR^2),\cE)$ to $(D([0,1],\IR^2),\cD)$ as follows: for any set 
$A\in \cD$, we have $A\cap C([0,1],\IR^2) =: B\in \cE$ and, via putting 
$\IW(A) = \IW(B)$ for $A\in \cD$, we retain $\IW(C([0,1],\IR^2))=1$. Also, on 
account of (1.15), for any map $h:D([0,1],\IR^2) \to \IR^2$ that is 
$\cD$ measurable, with $h$ as in (1.14), the $\IR^2$-valued map 
$h(\W(t)$diag$(\gamma^{-1/2},(1-\gamma^{-1})^{1/2}))$ is seen to be a two 
dimensional random vector. For further details concerning the above lines, we
refer to the middle three paragraphs of page 2948 of \Cs\ and Martsynyuk 
(2011). Also, in Section 1.2 of the latter paper, we say that a sequence of 
random elements $\{G_n(t),0\leq t\leq 1\}_{n\geq 1}$ of $(D([0,1],\IR^2),\cD)$ 
converges weakly to a continuous random element $\{G(t),0\leq t\leq 1\}$ of 
$(D([0,1],\IR^2),\cD)$, if, as $n\to\infty$, $g(G_n(t)) \dto g(G(t))$ for all 
functionals $g:D([0,1],\IR^2)\to\IR$ that are $\cD$ measurable and 
$\Delta$-continuous, or $\Delta$-continuous except on a set of measure zero 
in $\cD$ with respect to the measure generated by $G(t)$. We note that 
Corollary 1 also holds true in terms of the latter definition of weak 
convergence, i.e., with replacing $G_n(t)$ by $n^{-1/2} \Z_{[nt]}$, and $G(t)$ 
by $\W(t)$, in this definition in hand (cf.\ also Remarks 1.2, 1.3, and 
Proposition 1.1.\ of \Cs\ and Martsynyuk (2011)).

Back to Corollary 1, on taking $t=1$, it reduces to Theorem of Heyde (1993) 
(cf.\ Theorem B). 

As another immediate example of Corollary 1, as $n\to\infty$, under the 
conditions of Theorem 1, we conclude
$$
\left(n^{-1/2}\int^1_0 X_{[nt]}dt, \, 
n^{-1/2}\int^1_0 Y_{[nt]}dt \right)\phantom{XXXXXXXX}
$$

$$ 
\dto \,\left(\gamma^{-1/2} \int^1_0 W_1(t)dt, \, (1-\gamma^{-1})^{1/2} 
\int^1_0 W_2(t)dt \right)
$$
$$
 \deq \Big(\gamma^{-1/2}N_1(0,1/3),(1-\gamma^{-1})^{1/2} N_2(0,1/3)\Big),
$$
where $N_1(0,1/3)$ and $N_2(0,1/3)$ are independent normal random variables 
with mean 0 and variance 1/3. Here, and throughout, $\deq$ stands for 
equality in a distribution.

Define now the continuous time process $\Z_{nt} := \{X_{nt}, Y_{nt}\}$ by 
linear interpolation for $t\in [0,1]$, i.e., $\Z_{nt}$ are random elements 
of $C([0,1],\IR^2)$, and interpret the uniform $\Delta$-norm defined in
(1.10) for functions in this function space.

Recall the definition of the two dimensional Strassen (1964) class of 
absolutely continuous functions:
$$
\cS^{(2)} = \left\{\left({f}(x),{g}(x)\right),\, 0\leq x\leq 1: 
f(0)=g(0)=0,
\int^1_0 (\dot{f}^2(x) + \dot{g}^2(x))dx \leq 1\right\}.
$$

As noted already in the Introduction, via using an alternative construction of 
random walk $\{X_n,Y_n\}$ on ${\mathbb Z}^2$ with horizontal step transition 
probabilities $1/2-p_j$ and vertical step transition probabilities $p_j$, i.e.,
an inverted version of transition probabilities as in (1.1), and assuming 
condition (\ref{eq1.2}), Cs\'aki {\it et al.} (2013) established a joint 
strong approximation for $\{X_n,Y_n\}$ by two independent standard Wiener 
processes as in their Theorem 1.1. As a consequence of the latter Theorem 1.1 
and the two-dimensional version of the Strassen (1964) functional iterated 
logarithm law that is stated as Lemma D in Cs\'aki {\it et al.} (2013) 
for the continuous time process ${\mathbf Z}_{nt} :=\{X_{nt},Y_{nt}\}$ as 
right above, they concluded a functional law of the iterated logarithm as 
stated in their Corollary 1.1. Via our Theorem 2 parts (i)-(iii) of the latter 
corollary are inherited in the present context, and they read as follows.

\noindent{\bf Corollary 2} ~{\sl For the random walk $\Z_{n\cdot}$, as a
consequence of {\rm Theorem 2}, we conclude that
\begin{itemize}
\item[{\rm(i)}] the sequence of random vector-valued functions
$$
\left(\! \gamma^{1/2} \frac{X_{nx}}
 {(2n\log\log n)^{1/2}},
 \Bigg(\! \frac{\gamma}{\gamma\!-\!1}\!\Bigg)^{1/2}\!\!
\frac{Y_{nx}}
{(2n\log\log n)^{1/2}},   0\!\leq\! x\! \leq\! 1\! \right)_{\!n\geq 3}
$$
is almost surely relatively compact in $C([0,1],\IR^2)$ {\sl in the uniform
$\Delta$-norm topology}, and its limit points is the set of functions
$\cS^{(2)}$ (i.e., the collection of a.s.\ limits of convergent subsequences
{\sl in uniform $\Delta$-norm}).
\item[{\rm(ii)}] In particular, the vector sequence
$$
\left( \frac{X_n}{(2n\log\log n)^{1/2}}
\frac{Y_n}{(2n\log\log n)^{1/2}}\right)_{n\ge 3}
$$
is almost surely relatively compact in the rectangle
$$
\left[ -\frac{1}{\sqrt \gamma},\, \frac{1}{\sqrt \gamma}\right]
\times \left[ - \frac{\sqrt{\gamma-1}}{\sqrt \gamma}, \,
\frac{\sqrt{\gamma-1}}{\sqrt \gamma}\right]
$$
and the set of its limit points is the ellipse
$$
\Bigg\{(x,y): \gamma x^2 + \frac{\gamma}{\gamma-1}\, y^2\le 1\Bigg\}.
$$
\item[{\rm(iii)}] Moreover,
$$
\limsup_{n\to\infty} \frac{X_n}{\sqrt{2n\log\log n}} =
\frac{1}{\sqrt\gamma}\ \rm{a.s.}
$$
and
$$
 \limsup_{n\to\infty} \frac{Y_n}{\sqrt{2n\log\log n}} =
\frac{\sqrt{\gamma-1}}{\sqrt \gamma} \ \rm{a.s.}
$$
\end{itemize}}

For a survey on the path behaviour of random walks on a 2-dimensional 
anisotropic lattice with asymptotic density conditions {\sl \`a la} Heyde 
(1982, 1993) we refer to Cs\'aki {\sl et al.} (2015). 

\subsection{A nearest-neighbour random walk in a random environment on 
$\IZ^2$; Random column configurations}\label{ssect2}

Shuler (1979) formulated three conjectures on the asymptotic properties of a 
nearest-neighbour random walk on $\IZ^2$ that is allowed to make horizontal 
steps everywhere but vertical steps only on a random fraction of the columns. 

Following den Hollander (1994), let
\beq
C=\{C(x)\}_{x\in\IZ}
\label{eq1.16}
\eeq
be a random $\{0,1\}$-valued sequence with probability law $\mu$ on 
$\{0,1\}^\IZ$, satisfying the assumptions:
\begin{eqnarray}
(*)\ \ &&\mu \hbox{ is stationary and ergodic (w.r.t.\ translations in $\IZ$)},
\nonumber\\
&& 0 < q:= \mu(C(0)=1)\leq 1,
\label{eq1.17}
\end{eqnarray}
i.e.,  $q = E_\mu C(0)$, the expected value of $C(0)$ with respect to $\mu$.

We note in passing that due to the assumed stationarity of the probability 
measure $\mu$ on $\{0,1\}^\IZ$, we have (1.17) for $C(x)$ for all $x\in \IZ$ 
and $E_\mu C(x)=q$.

Given the sequence $C=\{C(x)\}_{x\in\IZ}$, a {\sl random environment} is 
constructed the following way: ~(i) the horizontal edges of the lattice 
$\IZ^2$ are left unchanged, i.e., all the horizontal edges in all the rows 
are kept, and (ii) all the vertical edges in the column $x$ are left or 
erased, depending on whether $C(x) = 1$ or $C(x)=0$. Thus, given 
$C=\{C(x)\}_{x\in\IZ}$, all the rows are connected, but only a part of the 
columns are.

Now, given $C=\{C(x)\}_{x\in\IZ}$, let
\beq
\left\{\Z(n)\right\}_{n\geq 0} = \{X(n),Y(n)\}_{n\geq 0}
\label{eq1.18}
\eeq
be the random walk that starts from the origin at time 0 and chooses with 
equal probability one of the unerased edges of the lattice $\IZ^2$ adjacent 
to the current position and jumps along it to a neighbouring site, i.e., 
given $C$, the random walk $\Z(n)$, $n\in \N$, on the path space 
$(\IZ^2)^{\N}$ with probability law $P_C$ as in (4) of den Hollander (1994) 
has the following transition probabilities:
\begin{eqnarray}
&& P_C(\Z(n+1)=(x\pm 1,y)|\Z(n)=(x,y))=1/2\quad \hbox{ if~ } C(x)=0\nonumber\\
&&P_C(\Z(n+1) = (x\pm 1,y\pm 1)|\Z(n)=(x,y))=1/4 \quad \hbox{ if~ } C(x) = 1.~~
\label{eq1.19}
\end{eqnarray}
On integrating over $C$ with respect to $\mu$, the thus obtained 
{\sl random walk in random environment} process has probability law 
$P := \int P_C\mu(dC)$.

We note in passing that with $q=1$, (1.19) becomes a simple symmetric random 
walk on $\IZ^2$. 

The three conjectures formulated by Shuler (1979), and referred to as 
``Ans\"atze'' by him, relate to the asymptotic behaviour of $\bZ(n)$ under 
the law $P$ as $n\to\infty$. They concern the total number of steps and the 
mean-square displacement in the $x$ (horizontal) and $y$ (vertical) 
directions, the probability of return to the origin, and the expected number 
of distinct sites visited, in the first $n$ steps of $\bZ(n)$.

For $n$ fixed, and keeping the notation used by den Hollander (1994), let
\begin{eqnarray}
n_x(n) &:=& \left|\{0\leq m<n:X(m+1) \neq X(m)\}\right|,\nonumber\\
n_y(n)  &:=& \left|\{0\leq m<n:Y(m+1) \neq Y(m)\}\right|,
\label{eq1.20}
\end{eqnarray}
respectively denote the total number of horizontal and vertical steps in the 
first $n$ steps of $\bZ(n)$.  Proved as {\sl Ansatz 1} of Shuler (1979) in 
den Hollander (1994), with $C=\{C(x)\}_{x\in\IZ}$ of (1.16) and $\mu$ as in 
(1.17), we have
\begin{eqnarray}
&& \lim_{n\to\infty} \frac{n_x(n)}{n} = q_x \quad P\hbox{--a.s.}\nonumber\\
&& \lim_{n\to\infty}\frac{n_y(n)}{n} = q_y\quad P\hbox{--a.s.}
\end{eqnarray}

Clearly, $n_x(n) + n_y(n)=n$ and, in view of (1.21), $q_x+q_y=1$. Also, as 
noted by den Hollander (1994), $q_x$ and $q_y$ denote the density of 
horizontal and vertical bonds in the lattice $\IZ^2$, i.e., $q_x+q_y=1$ and 
$q_y/q_x=q$ by $(*)$. Consequently, the respective limits in (1.21) are
\beq
q_x = 1/(1+q) \hbox{ ~and~ } q_y = q/(1+q).
\label{eq1.22}
\eeq

We note in passing that the $P$--a.s.\ conclusions of (1.21) can be viewed as 
respective analogues of those of (1.6) and (1.7) for the anisotropic random 
walk $\bZ_n$ under the asymptotic density condition (1.2) for its transition 
probabilities as in (1.1) in that, under $(*)$, the asymptotic behaviour of 
the random walk $\Z(n)$ with its $P_C$ transition probabilities as in (1.19) 
integrated over $C$ behaves like an anisotropic random walk on $\IZ^2$, with 
the anisotropy determined only by $q$ via $\mu$ of $(*)$. Moreover, $q$ 
viewed as the density of the connected columns via $q=q_y/q_x$, it 
``coincides'' with the so-called {\bf dimensional anisotropy} (cf.\ the 
quoted {\bf 3.~Final Remark} of Heyde (1993) right above Remark 1)
\beq
\lim_{n\to\infty} \frac{n^{-1}\ell(n)}{n^{-1}k(n)} = 
\frac{1-\gamma^{-1}}{\gamma^{-1}} = \gamma - 1
\label{eq1.23}
\eeq
of the anisotropic random walk $\bZ_n$ under the condition (1.2) that, in 
turn, implies the respective a.s.\ conclusions of (1.6) and (1.7) via the 
Proof of Theorem 1 of Heyde (1982), as noted right after (1.6).

\noindent{\bf Remark 4}\label{rem4}
~In view of our conclusion in Remark 1 concerning the so-called uniform 
periodic case, and our lines right above on $q=q_y/q_x$ {\sl versus} (1.23), 
it appears to be reasonable to say that $(*)$ together with (1.19), in the 
long run amount to a stationary ergodic randomization of deleting columns in 
the uniform periodic case on an average $q=E_\mu C(0) = EC(x)$ times, instead 
of $(\gamma-1)$ fraction of times as in (1.23) that, for a given positive 
integer $L\geq 1$, is equal to $1/L$ in the fixed column uniform periodic 
configuration case. On taking $q=1/L$ in $(*)$ in combination with the 
transition probabilities of (1.19) in the random column configuration model, 
in the long run amounts to realizing the dimensional anisotropy of a fixed 
column uniform periodic configuration $1/L$ times on an average, instead of 
having it equal to $(\gamma-1)=1/L$ eventually under the condition (1.2).

Further along these lines, from the $P$--a.s.\ limits of (1.21) that coincide 
with Ansatz 1 of den Hollander (1994), he concludes
\begin{eqnarray}
&&EX^2(n) \sim \frac{1}{1+q} n\nonumber\\ 
&&EY^2(n) \sim \frac{q}{1+q} n.
\label{eq1.24}
\end{eqnarray}

Here $E$ denotes expectation with respect to $P=\int P_C\mu(dC)$, and the 
statements themselves are analogous to (b) of Theorem B, established by Heyde 
(1993) under his condition (1.2) with $\eta = 0$.

Under the additional {\sl mixing assumption} (19) of den Hollander (1994) on 
$\mu$, he concludes
\beq
Q(n) := P\left( \bZ(2n) = (0,0)\right) \sim \frac{1+q}{2\pi q^{1/2} n}
\label{eq1.25}
\eeq
and
\beq
R(n) := E|\{\bZ(0), \bZ(1),\ldots,\bZ(n)\}| \sim 
\frac{2\pi q^{1/2}  n}{(1+q)\log n},
\label{eq1.26}
\eeq
where $E$ denotes expectation with respect to $P$ as in (1.24), and the 
integrand $|\{\bZ(0), \bZ(1),\ldots,\bZ(n)\}|$ stands for the number of 
distinct sites visited, as in den Hollander (1994).

As noted before, with $q=1$, (1.19) becomes a simple symmetric random walk on 
$\IZ^2$, and (1.25) and (1.26) are well known results in this classical context.

Also, as noted already, the conclusions of (1.21) constitute {\sl Ansatz 1} 
of Shuler (1979). We also note in passing that the respective results of 
(1.25) and (1.26) that are based on those of (1.24), that is (10) in den 
Hollander (1994), coincide with (11) and (12) of den Hollander (1994), which 
in turn  respectively amount to {\sl Ansatz 2} and {\sl Ansatz 3} of Shuler 
(1979). Concerning the latter (10), (11) and (12) of his paper, den Hollander 
(1994) writes as follows:

\begin{quote}

``Note that the coefficients in (10)-(12) depend on $\mu$ only via the density 
of connected columns $q$.  Moreover, (10)-(12) are precisely what one would 
find for a random walk on the full lattice that makes horizontal and vertical 
steps with probability $q_x = 1/(1+q)$, resp. $q_y=q/(1+q)$. Actually, this 
was the {\sl main idea behind the formulation of the Ans\"atze} in the first 
place: asymptotically the random walk should behave like an anisotropic 
random walk on $\IZ^2$, with the anisotropy determined {\sl only by $q$ and 
not by any other parameters of $\mu$}. Though, as we shall see later, this 
is {\sl not} quite true in full generality under $(*)$ [not for $Q(n)$ and 
$R(n)$ at least], it is indeed true for many asymptotic properties associated 
with the random walk and within a large class of column distributions 
(including periodic distributions and i.i.d.\ distributions)."

\end{quote}

The in here  mentioned condition $(*)$ of den Hollander (1994) is that of 
(1.17) on $\mu$ above, and the hinted at further condition is the 
{\sl mixing assumption} of (19) on $\mu$ in den Hollander (1994) that we 
already mentioned fleetingly when introducing (1.25) and (1.26) above.

\noindent{\bf Observation}
In this section we are to spell out analogues of Theorems 1 and 2 under the 
$(*)$ den Hollander (1994) condition of (1.17) on $\mu$.  These analogues 
will constitute evidence to saying that, under the latter condition, 
asymptotically the random walk $\bZ(n)$ in random environment with probability 
law $P=\int P_C\mu(dC)$ (cf.\ (1.16), (1.19) respectively for $C$ and $P_C$) 
behaves like an anisotropic random walk on $\IZ^2$ does, with the anisotropy 
determined only by $q$ and not by any other parameters of $\mu$, as if it 
were an anisotropic 2-dimensional random walk under the asymptotic density 
condition of Heyde (1982), with $\gamma = 1+q$, $0<q\leq 1$, in (1.2) for the 
transition probabilities of (1.1).

As in den Hollander (1994), let $B_q=\{B_q(t)\}_{t\geq 0}$ be anisotropic 
Brownian motion on $\IR^2$ with diffusion matrix
\beq
D = \left( \begin{array}{cc}
1/(1+q) & 0\\
0 & q/(1+q) \end{array} \right), 
\label{eq1.27}
\eeq
i.e., $EB_q(t) = (0,0)$ and the finite dimensional distributions of 
$\{B_q(t)\}_{t\geq 0}$ are Gaussian with covariance matrix 
min$(s,t)D, ~0\leq s,t < \infty$. Consequently, $\{B_q(t)\}_{t\geq 0}$ has 
the following representation
\begin{eqnarray}
\{B_q(t)\}_{t\geq 0}& \deq& 
\left\{\bW(t)\diag\left((1/(1+q))^{1/2}, 
(q/(1+q))^{1/2}\right)\right\}_{t\geq 0}\nonumber\\
&\deq&\left\{W_1(t/(1+q)),W_2(tq/(1+q))\right\}_{t\geq 0}
\label{eq1.28}
\end{eqnarray}
with $\{\bW(t)\}_{t\geq 0} := \{(W_1(t),W_2(t))\}_{t\geq 0}$, where 
$\{W_1(t), \, t\geq 0\}$ and $\{W_2(t),\, t\geq 0\}$ are two independent 
standard Wiener processes. 

We now reformulate the invariance principle of Theorem 1 of den Hollander 
(1994) for $\bZ(n)$ {\sl \`a la} Theorem 1 in our previous section.

\noindent {\bf Theorem C} {\sl Assume $(*)$ for the probability 
law $\mu$ of the random $\{0,1\}$-valued sequence $\{C(x)\}_{x\in\IZ}$ on 
$\{0,1\}^\IZ$. Then $(1.21)$ is true and, on an appropriate probability 
space for $\bZ(n)$ on $\IZ^2$ with transition probabilities as in $(1.19)$ 
and probability law $P:=\int P_C\mu (dC)$, one can construct two independent 
standard Wiener processes $\{W_1(t),\, t\geq 0\}$, $\{W_2(t),\, t\geq 0\}$ so 
that, as $n\to\infty$, with}
\beq
\{\bW_n(t),t\geq 0\}_{n\geq 0} := 
\left\{ \frac{W_1(nt/(1+q))}{n^{1/2}}, 
\frac{W_2(ntq/(1+q))}{n^{1/2}}, \, t\geq 0\right\}_{n\geq 0}
\label{eq1.29}
\eeq
{\sl and} $q$ {\sl as in} (1.17), {\sl we have for all fixed} $T>0$
\begin{eqnarray}
\lefteqn{\sup_{0\leq t\leq T} \Vert n^{-1/2}\bZ([nt])-\bW_n(t)\Vert}\nonumber\\
&=& \sup_{0\leq t\leq T} \left\Vert
\frac{X([nt])-W_1(nt/(1+q))}{n^{1/2}},
\frac{Y([nt])-W_2(ntq/(1+q))}{n^{1/2}} \right\Vert\nonumber\\
&=& o_P(1).
\label{eq1.30}
\end{eqnarray}

Moreover, {\sl under the same condition}, i.e., only under the $(*)$
den Hollander (1994) condition of (1.17) on $\mu$, we also have the 
following Strassen type extended version of (21) of his Theorem 2 
({\sl ibid.}).

\noindent{\bf Theorem 3} {\sl Assume $(*)$ as in {\rm(1.17)}. Then
on the probability space of {\rm Theorem C} with the same independent standard 
Wiener processes $W_1$ and $W_2$, we have}
\begin{eqnarray}
\lefteqn{\sup_{0\leq t\leq 1}\left\Vert (2n\log\log n)^{-1/2} \bZ([nt])
-(2\log\log n)^{-1/2}\bW_n(t)\right\Vert}\nonumber\\
&=& \sup_{0\leq t\leq 1}
\left\Vert\frac{X([nt])-W_1(nt/(1+q))}{(2n\log\log n)^{1/2}},
\frac{Y([nt])-W_2(ntq/(1+q))}{(2n\log\log n)^{1/2}} \right\Vert 
\nonumber\\[2ex]
&=& o(1) \quad P\hbox{--a.s.}, \, \,  n\to\infty.
\label{eq1.32}
\end{eqnarray}

\noindent{\bf Remark 5 } ~Under their respective conditions, 
Theorem 1 and Theorem C  ``coincide", and so do also Theorem 2 and Theorem 3. 
Thus, to the extent of simultaneously having a weak Donsker and a strong 
Strassen type asymptotic behaviour, the two random walks in hand behave 
similarly. Moreover, just like having Theorem 1 under the weaker conditions 
of (1.8) than that of (1.6), Theorem C will also be seen to be true under 
assuming only convergence in probability versions of the $P$--a.s. conclusions 
of (1.21).  We do not however know how to go about weakening the $(*)$ 
assumption of (1.17) for the probability law $\mu$ so that it would only 
yield the desired weaker version of (1.21).

\noindent{\bf Remark 6} ~Mutatis mutandis, Corollaries 1 and 2 also 
hold true in the context of Theorems C and 3 respectively. 

\section{Proofs}\label{sect2}

\setcounter{equation}{0}

\subsection{Preliminaries and proofs of Theorems 1 and 2}\label{ssec2.1}

As right after Corollary A {\it \`a la} 3.~Proofs of Heyde (1982), let again
$\sigma_0 = 0 < \sigma_1 < \sigma_2 < \cdots$ be the successive times at 
which the values of the $X_i - X_{i-1}$, $i=1,2,\ldots$ are nonzero, and put 
again  $S_1(k) = X_{\sigma_k}$.  By the assumed symmetry of the transition 
probabilities in (1.1), $\{S_1(k),\, k\geq 0\}$ is a simple symmetric random 
walk on $\IZ$.  Also, $X_n = X_{\sigma_k}$ for 
$\sigma_k \leq n < \sigma_{k+1}$. For $n$ fixed let
$$
\sigma_{k(n)} := \hbox{max}[j:j\leq n, \, X_j \neq X_{j-1}].
$$
Then again, as in (1.5),
\beq
X_n = X_{\sigma_{k(n)}} = S_1(k(n))
\label{eq2.1}
\eeq
is the horizontal position of the walk $\bZ_n = (X_n,Y_n)$ after $k(n)$ 
horizontal steps in the first $n$ steps of $\bZ_n$.

Clearly, in view of (2.1), $\ell(n) := n-k(n)$ is the number of vertical steps 
in the first $n$ steps of $\bZ_n$. On its own, the vertical position $Y_n$ 
of the walk $\bZ_n = (X_n,Y_n)$ can be dealt with similarly to that of its 
horizontal position $X_n$.  As in the proof of Theorem in Heyde (1993), let 
$\tau_0 = 0 < \tau_1<\tau_2<\cdots$ be the successive times at which the 
values of $Y_i - Y_{i-1}$, $i=1,2,\ldots$ are nonzero and put 
$S_2(k)=Y_{\tau_k}$. Then, again by the assumed symmetry of the transition 
probabilities (1.1), $\{S_2(k),\, k\geq 0\}$ is a simple symmetric random 
walk on $\IZ$, and $Y_n=Y_{\tau_k}$ for $\tau_k\leq n<\tau_{k+1}$. For $n$ 
fixed, let
$$
\tau_{\ell(n)} := \max[j:j\leq n,\, Y_j\neq Y_{j-1}].
$$
Then, for $Y_n$, the vertical position of the random walk $\bZ_n = (X_n,Y_n)$ 
after $n$ steps, we have
\beq
Y_n=Y_{\tau_{\ell(n)}} = S_2(\ell(n))
\label{eq2.2}
\eeq
after $\ell(n)$ vertical steps in the first $n$ steps of $\bZ_n$.

As above, with $\sigma_k$ and $\tau_k$  standing for the times at which $X$ 
and $Y$ make their respective $k$-th steps $(\sigma_0 = \tau_0 = 0)$, both
\beq
\{X_{\sigma_k}, k\geq 0\} = \{S_1(k), k\geq 0\}
\label{eq2.3}
\eeq
\beq
\{Y_{\tau_k}, k\geq 0\} = \{S_2(k), k\geq 0\}
\label{eq2.4}
\eeq
are simple symmetric random walks on $\IZ$, and are 
{\sl defined independently of each other}. Consequently, the random walk 
$\bZ_n = (X_n,Y_n)_{n\geq 0}$ with transition probabilities as in (1.1) and 
viewed {\it \`a la} (2.1) and (2.2) can be studied in terms of 
\beq
\bZ_n = (X_n,Y_n) =(X_{\sigma_{k(n)}},Y_{\tau_{\ell(n)}}) = 
(S_1(k(n)),S_2(\ell(n))
\label{eq2.5}
\eeq
where, as before, $k(n)$ and $\ell(n)$ are the respective numbers of 
horizontal and vertical steps in the first $n$ steps of $\bZ_n$, 
$S_1(\cdot)$ and $S_2(\cdot)$ are independent simple symmetric random walks 
on $\IZ$, and $k(n) + \ell(n) = n$. 

Now, without changing their distribution, on an appropriate probability space 
for the two independent simple symmetric random walks 
$\{S_i(j),j\geq 0\}$, $i=1,2$, one can construct two independent standard 
Wiener processes $\{W_1(t),t\geq 0\}$ and $\{W_2(t),t\ge 0\}$ so that 
(cf.\ Koml\'os, Major and Tusn\'ady [KMT] (1975))
\beq
|S_i(j)-W_i(j)| = O(\log j) \qas,~ i=1,2,
\label{eq2.6}
\eeq
as $j\to\infty$, and
\beq
\sup_{0\le t\le 1} |S_i([nt]) - W_i(nt)| = O(\log n) \qas,~ i=1,2,
\label{eq2.7}
\eeq
as $n\to\infty$.

Consequently, as $\noo$, the components of $\bZ_n$ as in (2.5) can be studied 
independently via the approximations
\beq
|S_1(k(n))-W_1(k(n))| = O(\log k(n)) \qas,
\label{eq2.8}
\eeq
\beq
\sup_{0\le t\le 1}|S_1(k(nt))-W_1(k(nt))| = O(\log k(n)) \qas,
\label{eq2.9}
\eeq
and
\beq
|S_2(\ell(n))-W_2(\ell(n))| = O(\log \ell(n)) \qas,
\label{eq2.10}
\eeq
\beq
\sup_{0\le t\le 1}|S_2(\ell(nt))-W_2(\ell(nt))| = O(\log \ell(n)) \qas,
\label{eq2.11}
\eeq
where $k(nt)$ and $\ell(nt)$ are the respective numbers of horizontal and 
vertical steps in the first $[nt]$ steps of $\bZ_{[nt]}$, $0\le t\le 1$, and 
$k(nt)+\ell(nt) = [nt]$. 

Now, under condition (1.2) we have (1.6) and, as in the Proof of Theorem 1 of 
Heyde (1982), on writing $k(n) = n\gamma^{-1}(1+\varepsilon_n)$, with 
$\ep_n\to 0$ a.s.\ as $n\to\infty$, we arrive at (cf.\ (2.8))
\begin{eqnarray}
|S_1(k(n)) - W_1(n\gamma^{-1}(1+\ep_n))| &&=~ 
O(\log n\gamma^{-1}(1+\ep_n)) \qas\nonumber \\
&& = ~O(\log n) \qas
\label{2.12}
\end{eqnarray}
as $n\to\infty$. Moreover, as $n\to\infty$, we have also (cf.\ (2.9))
\beq
\sup_{0\leq t\leq 1}\Big|S_1(k(nt)) - W_1(nt \gamma^{-1}(1+\ep_{n}))\Big|
= O(\log n) \qas
\label{(2.12)}
\eeq

Recalling that for $X_n$, the horizontal position of the random walk after 
$n$ steps, we have (cf. (2.1)) $X_n = X_{\sigma_{k(n)}} = S_1(k(n))$, via 
(2.12), as $n\to\infty$, we arrive at  
$$
\gamma^{1/2}X_n = W_1(n(1+\ep_n)) + O(\log n) \qas
$$
that in turn rhymes with (1.3) of Theorem A, a result of Heyde (1982), except 
the a.s. $O(\cdot)$ rate, that in his case is due to using a result of 
Strassen (1967), that amounts to the best possible approximation of the 
partial sums $S_1(\cdot)$ by Brownian motion via Skorokhod stopping times, 
while that of KMT (1975) that we use amounts to its globally best possible 
approximation.  As we will however see later on, both serve equally well in 
the present context.

In view of (2.13), as $n\to\infty$, we have
\begin{equation}
~~~~\sup_{0\leq t\leq 1} \left| \frac{X_{[nt]}}{n^{1/2}} - 
\frac{W_1(nt\gamma^{-1}(1+\ep_{n}))}{n^{1/2}} \right|
= O\left(\frac{\log n}{n^{1/2}} \right) \hbox{~ a.s.}
\label{eq2.14}
\end{equation}

As to $Y_n$, the vertical position of the random walk $\bZ_n=(X_n,Y_n)$, 
recall (cf.\ (2.5)) that we have $Y_n=Y_{\tau_{\ell(n)}} = S_2(\ell(n))$, 
as well as $\bZ_n=(X_n,Y_n) =$ $ (S_1(k(n)),S_2(\ell(n)))$, and 
$k(n) + \ell(n) = n$, the sum of the number of horizontal and vertical steps 
in the first $n$ steps of $\bZ_n$.  Hence, in view of (1.7), and on writing 
$\ell(n) = n(1-\gamma^{-1})(1+\tilde\ep_n)$, with $\tilde\ep_n\to 0$ a.s.  
as $n\to\infty$, {\sl mutatis mutandis} in concluding (2.14), with a standard 
Wiener process $\{W_2(t), t\!\geq\! 0\}$ as in (2.6), we arrive at
\begin{equation}
\sup_{0\leq t\leq 1}\!\left| \frac{Y_{[nt]}}{n^{1/2}} - 
\frac{W_2(nt(1\!-\!\gamma^{-1})(1\!+\!\tilde\ep_{n}))}{n^{1/2}} \right| = 
O\left(\!\frac{\log n}{n^{1/2}}\! \right) \qas
\label{eq2.15}
\end{equation}
as $n\to\infty$.

We note again that the conclusions of (2.12)--(2.15) are based on the Heyde 
(1982) condition (1.2) yielding (1.6) and (1.7), and, consequently, on writing
\begin{equation}
k(n) = n\gamma^{-1}(1+\ep_n), ~\ell(n) = n(1-\gamma^{-1})(1+\tilde\ep_n)
\label{2.16}
\end{equation}
with $\ep_n$ and $\tilde\ep_n$ both converging almost surely to 0 as 
$n\to\infty$.

We now recall that under the weaker condition that $\eta = 0$ in (1.2), 
Heyde (1993) concluded (1.6) and (1.7) to hold true in probability as 
indicated in (1.8). Consequently, we can again write $k(n)$ and $\ell(n)$ as 
above in (2.16), but now with $\ep_n$ and $\tilde \ep_n$ both converging in 
probability to 0.  This, in turn, results in having (2.12)--(2.15) in 
probability instead of almost surely. In particular, for further use, we spell 
out the in probability versions of (2.14) and (2.15) as follows.

With the two independent standard Wiener processes $W_1$ and $W_2$ as in 
(2.6), and writing $k(n)$ and $\ell(n)$ as in (2.16) with $\ep_n$ and 
$\tilde\ep_n$ both converging in probability to 0 as $n\to\infty$ on account 
of having (1.8) under the condition that $\eta = 0$ in (1.2), as $n\to\infty$, 
we arrive at
\begin{equation}
\sup_{0\leq t\leq 1} \left| \frac{X_{[nt]}}{n^{1/2}} - 
\frac{W_1(nt\gamma^{-1}(1+\ep_{n}))}{n^{1/2}} \right|
= o_P(1)
\label{2.17}
\end{equation}
and
\begin{equation}
\sup_{0\leq t\leq 1}\left| \frac{Y_{[nt]}}{n^{1/2}} - 
\frac{W_2(nt(1\!-\!\gamma^{-1})(1\!+\!\tilde\ep_{n}))}{n^{1/2}} \right| 
= o_P(1).
\label{2.18}
\end{equation}

\noindent{\bf Lemma 1} {\sl Under the condition $(1.2)$ as is, that 
via $(2.16)$ yields the a.s.\ convergence of both $\ep_n$ and $\tilde\ep_n$ 
to $0$ as $n\to\infty$, and also  under the condition $(1.2)$ with 
$\eta = 0$, that via $(2.16)$ results in $\ep_n$ and $\tilde\ep_n$ both 
converging in probability to 0 as $n\to\infty$, we have}
\begin{equation}
\sup_{0\leq t\leq 1} \left| \frac{W_1(nt\gamma^{-1}(1+\ep_{n}))-
W_1(nt\gamma^{-1})} {n^{1/2}} \right| = o_P(1)
\label{2.19}
\end{equation}
{\sl as well as}
\begin{eqnarray}
\sup_{0\leq t\leq 1} \left| 
\frac{W_2(nt(1\!-\!\gamma^{-1})(1\!+\! \tilde\ep_{n}))-
W_2(nt(1-\gamma^{-1}))} {n^{1/2}} \right|
=o_P(1).
\label{2.20}
\end{eqnarray}

We note in passing that these two statements hold true in terms of any 
generic standard Wiener process $W$. The present forms are only for the sake 
of easy reference when  proving Theorem 1, and, later on, when indicating the 
proof of Theorem C.

\noindent{\bf Proof of Lemma 1} First assume (1.2).  
Then with $\ep_n \to 0$ a.s.\  as $n\to\infty$, we have
\begin{eqnarray}
&&\sup_{0\leq t\leq 1} \left| \frac{W_1(nt\gamma^{-1}(1+\ep_{n}))-
W_1(nt\gamma^{-1})} {n^{1/2}} \right|\nonumber\\
&& \deq \sup_{0\leq t\leq 1} 
\left| W_1(t\gamma^{-1}(1+\ep_{n}))-W_1(t\gamma^{-1}) \right|, 
\hbox{ ~for each~ } n\geq 1,\nonumber\\
&& \deq \sup_{0\leq t\leq 1} 
\left| \gamma^{-1/2}(W_1(t+t\ep_{n})-W_1(t))\right|,  
\hbox{ ~for~ } 1< \gamma < \infty,  \nonumber\\
&& \leq \sup_{0\leq t\leq 1}\sup_{0<s\leq |\ep_n|} \gamma^{-1/2} 
\left| W_1(t+s)-W_1(t) \right|\nonumber\\
&& \leq \sup_{0\leq t\leq 1}\sup_{0<s\leq \ep} \gamma^{-1/2} 
\left| W_1(t+s)\!-\!W_1(t) \right| \qas \nonumber\\
&&= O(1)(\ep \log 1/\ep)^{1/2} \qas\nonumber\\
&&= o(1) \qas
\label{2.21}
\end{eqnarray}
by the L\'evy modulus of continuity as $\ep\downarrow 0$, and in view
of having the last inequality with any $\ep\!>\!0$, however small, for all 
but a finite number of $n$ on account of $\ep_n \to 0$ a.s. as $n\to\infty$.
Thus, via (2.21), we conclude (2.19)  under the condition (1.2), and a 
similar argument yields (2.20) as well under the same condition.

Assuming (1.2) with $\eta=0$, results in having $\ep_n\to 0$ in probability, 
as $n\to\infty$.  Consequently, and equivalently, every subsequence 
$\{\ep_{n_k}\}$ contains a further subsequence,
say $\{\ep_{n_{k(m)}}\}$, so that $\ep_{n_{k(m)}}\to 0$ a.s.\ as 
$n_{k(m)} \to \infty$.

Now to conclude (2.19) in this case, as in the proof above, we arrive at the 
first inequality with $\ep_n\to 0$ in probability, as $n\to\infty$, in the 
expression
\begin{equation}
\sup_{0\leq t\leq 1}\sup_{0<s\le |\ep_n|} \gamma^{-1/2}|W_1(t+s)-W_1(t)| =: 
\delta_n,
\label{2.22}
\end{equation}
and, in order to have (2.19), we wish to conclude that, as $n\to\infty$, 
$\delta_n\to 0$ in probability.  This, in turn, is seen to be equivalent to 
showing that $\delta_{n_{k(m)}}\to 0$ a.s. as $n_{k(m)}\to\infty$, on account 
of having $\ep_n \to 0$ in probability as $n\to\infty$, if and only if, as 
above, $\ep_{n_{k(m)}}\to 0$ a.s. as $n_{k(m)} \to \infty$.

Accordingly, we consider
\begin{eqnarray}
\delta_{n_{k(m)}}& :=& 
\sup_{0\le t\le 1} \sup_{0<s<|\ep_{n_{k(m)}}|} 
\gamma^{-1/2}|W_1(t+s)-W_1(t)|\nonumber\\
&\le &  \sup_{0\le t\le 1} \sup_{0<s<\ep} \gamma^{-1/2}|W_1(t+s)-W_1(t)| \qas, 
\nonumber\\
&=& O(1)(\ep \log 1/\ep)^{1/2} \qas\nonumber\\
&=& o(1) \qas
\label{2.23}
\end{eqnarray}
by the L\'evy modulus of continuity as $\ep\downarrow 0$, and on account of 
having the last inequality with any $\ep>0$, however small, for all but a 
finite number of the sequence $\{n_{k(m)}\}$ that increases to infinity 
and $\ep_{n_{k(m)}} \to 0$. This amounts to having $\delta_n$ of (2.22) 
converging to 0 in probability, as $n\to\infty$, that in turn yields (2.19) 
under the condition (1.2) with $\eta=0$. A similar argument results in having 
(2.20) as well under the same condition. This also completes the proof of 
Lemma 1. $\quad\Box$

\noindent{\bf Proof of Theorem 1}  With $W_1$ and $W_2$ as in (2.6), 
and assuming (1.2) with $\eta = 0$, we combine (2.17) with (2.19) and (2.18) 
with (2.20), and thus conclude Theorem 1 with $T=1$, without loss of 
generality, i.e., similarly for all fixed $T>0$ as well. $\quad\Box$

\noindent{\bf Lemma 2} {\sl Under condition $(1.2)$, as $n\to \infty$, 
we have}
\begin{equation}
\sup_{0\le t\le 1} \left| 
\frac{W_1(nt\gamma^{-1}(1+\ep_{n}))-W_1(nt\gamma^{-1})} 
{(2n\log\log n)^{1/2}} \right| = o(1) \qas
\label{2.24}
\end{equation}
{\sl as well as}
\begin{equation}
\sup_{0\leq t\leq 1} 
\left| \frac{W_2(nt(1\!-\!\gamma^{-1})(1\!+\! \tilde\ep_{n}))
-W_2(nt(1-\gamma^{-1}))} {(2n\log\log n)^{1/2}} \right|=o(1) \qas
\label{2.25}
\end{equation}
{\sl where $\ep_n$ and $\tilde\ep_n$ both converge almost surely to $0$, as 
indicated in $(2.16)$.}

\medskip
Just like in the case of Lemma 1, here too, $W_1$ and $W_2$ can be any 
generic Wiener processes.  The specific forms of (2.24) and (2.25) are only 
for the sake of convenient reference when proving Theorem 2 and, later on, 
when indicating the proof of Theorem 3.

In the proof of Lemma 2, we make use of the following large increment result 
of Cs\"org\H{o} and R\'ev\'esz (1979, 1981) for a standard Wiener process 
$\{W(t),t\!\geq\! 0\}$: {\sl  Let $a_T$ ($T\!\geq\! 0$) be a monotonically 
nondecreasing function of $T$ so that $T/a_T$ is nondecreasing and 
$0<a_T\leq T$.  Define $\beta_T=\{2a_T(\log \frac{T}{a_T} + 
\log\log T)\}^{-1/2}$.  Then}
\begin{equation}
~~~\limsup_{T\to\infty} \sup_{0\leq t\leq T-a_T} 
\sup_{0\leq s\leq a_T} \beta_T|W(t+s)-W(t)|=1 \qas
\label{eq2.26}
\end{equation}

\noindent
{\bf Proof of Lemma 2}  We have
\begin{eqnarray}
&&\sup_{0\leq nt\gamma^{-1} \leq n\gamma^{-1}} 
|W_1(nt\gamma^{-1}(1+\ep_{n})) - W_1(nt\gamma^{-1})|\nonumber\\
&& \leq \sup_{0\leq nt\leq n} \sup_{0\leq nt|\ep_{n}|\leq n |\ep_n|}
|W_1(nt\gamma^{-1}+nt\ep_{n}\gamma^{-1})-W_1(nt\gamma^{-1})|\nonumber\\
&& \leq \sup_{0\leq nt\leq n} \sup_{0\leq s\leq n\ep} 
|W_1(nt\gamma^{-1}+s\gamma^{-1})-W_1(nt\gamma^{-1})|
\label{eq2.27}
\end{eqnarray}
with any $\ep > 0$, however small, for all but a finite number of $n$, on 
account of $\ep_n \to 0$ a.s. as $n\to\infty$.

Moreover, with any $0<\ep < 1$,
\begin{eqnarray}
&&  \sup_{0\leq nt\leq n} ~\sup_{0\leq s\leq n\ep}  
|W_1(nt\gamma^{-1}+s\gamma^{-1})-W_1(nt\gamma^{-1})|\nonumber\\
&& \leq\sup_{0\leq nt\leq n-n\ep}\, 
\sup_{0\leq s\leq n\ep}|W_1(nt\gamma^{-1}+s\gamma^{-1})-W_1(nt\gamma^{-1})|
\nonumber\\
&& \quad  + \sup_{n-n\ep\leq nt\leq n}\, 
\sup_{0\leq s\leq n\ep}|W_1(nt\gamma^{-1}+s\gamma^{-1})-W_1(nt\gamma^{-1})|
\nonumber\\
&& = O\left(\left(2n\ep\left(\log \frac{n}{n\ep} + 
\log\log n\right)\right)^{1/2}\right)\nonumber\\
&&= O\left(\left(2n\ep\log \frac{1}{\ep} +
\ep 2n \log\log n\right)^{1/2}\right)\nonumber\\
&&= \Bigg(\ep\log \frac{1}{\ep} + \ep\Bigg)^{1/2}
O\left( (2n\log\log n)^{1/2}\right) \hbox{ a.s.}
\label{eq2.28}
\end{eqnarray}
as $n\to\infty$, on applying (2.26) twice.

On combining now (2.28) with (2.27) via letting $n\to\infty$ and 
$\ep \downarrow 0$, we arrive at
$$
\supt\frac{|W_1(nt\gamma^{-1}(1+\ep_{n}))-W_1(nt\gamma^{-1})|}
{(2n\log\log n)^{1/2}} = o(1) \qas,
$$
i.e., we have (2.24).  Mutatis mutandis, the conclusion of (2.25) is seen to 
be true as well along similar lines. \quad $\Box$

\noindent{\bf Proof of Theorem 2} ~By  (2.14), as $n\to\infty$,
\begin{eqnarray}
 \supt \left| \frac{X_{[nt]} - W_1(nt\gamma^{-1}(1+\ep_{n}))}
{(2n\log\log n)^{1/2}} \right|
= O\left( \frac{\log n}{(2n\log\log n)^{1/2}} \right) \hbox{ a.s.}
\label{eq2.29}
\end{eqnarray}
Putting (2.24) and (2.29) together, as $n\to\infty$, we conclude
\begin{eqnarray}
&& \supt \left| \frac{X_{[nt]} - W_1(nt\gamma^{-1})}
{(2n\log\log n)^{1/2}} \right|
=o(1)  \hbox{ a.s.}
\label{eq2.30}
\end{eqnarray}

In a similar fashion, as $n\to\infty$, one concludes
\begin{eqnarray}
&& \supt \left| \frac{Y_{[nt]} - W_2(nt(1-\gamma^{-1}))}
{(2n\log\log n)^{1/2}} \right|
=o(1)  \hbox{ a.s.}
\label{eq2.31}
\end{eqnarray}
as well. Consequently, via (2.30) and (2.31), we arrive at having 
Theorem 2. \quad $\Box$

\subsection{Preliminaries and proofs of Theorem C and Theorem 3}

Given $C=\{C(x)\}_{x\in\IZ}$ as in (1.16) with probability law $\mu$ on
$\{0,1\}^{\IZ}$ as in (1.17), i.e., assuming the $(*)$ condition of (2) of 
den Hollander (1994), both components of the random walk 
$\{\bZ(n)\}_{n\geq 0}=\{X(n),Y(n)\}_{n\geq 0}$ as in (1.18) with transition
probabilities as in (1.19) behave as simple random walks on $\IZ$, except for 
random time delays in the respective total number of horizontal and
vertical steps $n_x(n)$ and $n_y(n)$ in the first $n$ steps of $\bZ(n)$ as 
they are defined in (1.20). In this regard we quote from the first paragraph 
of Section 2 of den Hollander ({\sl ibid.}):

\begin{quote}

"Indeed, first $X(n)$ makes a succession of steps until it hits a connected
column. Next the walk spends some time on this column, during which $X(n)$
remains fixed and $Y(n)$ makes a succession of steps, until the walk decides 
to move off the column. Then $X(n)$ again takes over, until it hits a next
connected column, etc." 

\end{quote} 

Thereby, along the lines of the proof of {\sl Ansatz 1} ({\sl ibid.}),
the random walk $\{\bZ(n)\}_{n\geq 0}$ as in (1.18) with transition
probabilities as in (1.19) can be studied in terms of
\begin{eqnarray}
\{\bZ(n)\}_{n\geq 0} &=&\{X(n),Y(n)\}_{n\geq 0}\nonumber\\
                    &=&\{S_x(n_x(n)), S_y(n_y(n))\}_{n\geq 0},
\label{eq2.32)} 
\end{eqnarray}
where $S_x$ and $S_y$ are independent simple random walks, and $n_x(n)$
and $n_y(n)$ respectively stand for the total number of horizontal and 
vertical steps in the first $n$ steps of $\bZ(n)$ as defined in (1.20). 

Given $C=\{C(x)\}_{x\in\IZ}$ as in (1.16) with probability law $\mu$ on
$\{0,1\}^{\IZ}$ as in (1.17), i.e. assuming the $(*)$ condition of (2)
of den Hollander (1994), we have (1.21), proved ({\sl ibid.}) as 
{\sl Ansatz 1} of Shuler (1979), with the conclusions of (1.22). Hence, we 
may put 
\beq
n_x(n)=n(1/(1+q))(1+\varepsilon_n)
\label{eq2.33}
\eeq
and, consequently,
\beq
n_y(n)=n(q/(1+q))(1+\varepsilon_n)
\label{eq2.34}
\eeq
with $\varepsilon_n\to 0$ $P$-a.s., as $n\to\infty$.

As a consequence of having (2.32) in combination with (2.33) and (2.34),
one can construct two independent standard Wiener processes 
$\{W_1(t),\, t\geq 0\}$ and $\{W_2(t),\, t\geq 0\}$ so that via [KMT]
(1975), we arrive at
$$
\sup_{0\leq t\leq T} \left\Vert
\frac{X([nt])-W_1(nt(1/(1+q))(1+\varepsilon_n))}{n^{1/2}},
\frac{Y([nt])-W_2(nt(q/(1+q))(1+\varepsilon_n))}{n^{1/2}} \right\Vert
$$
\beq
= O(\log n/n^{1/2}) \, \, P\hbox{--a.s.} 
\label{eq2.35}
\eeq
for all fixed $T>0$, as $n\to\infty$.

Also, under the same conditions that yield (2.33) and (2.34), mutatis mutandis
accordingly, the respective statements of Lemmas 1 and 2 of Section 2.1 
continue to hold true in the present context.

Thus, on combining the latter version of Lemma 1 with the conclusion of
(2.35) with $T=1$, we obtain Theorem C with $T=1$,
without loss of generality, i.e., similarly for all fixed
$T>0$ as well (cf. also Remark 5). \quad $\Box$

Furthermore, the conclusion of Lemma 2 in the present context in
combination with that of (2.35) with $T=1$, yields Theorem 3.
\quad $\Box$

\bigskip\noindent
{\bf Acknowledgements} We wish to thank a number of referees for carefully
reading our manuscript, and for their constructive remarks and suggestions
that have very much helped us in preparing this revised version of our paper.
The research of E. Cs\'aki and P. R\'ev\'esz was supported by Hungarian 
National Research, Development and Innovation Office - NKFIH K 108615. 
The research of M. Cs\"org\H o was supported by an NSERC Canada Discovery 
Grant at Carleton University. The research of A. F\"oldes was supported by 
PSC CUNY Grant, No. 69040-0047.


\begin{thebibliography}{9}

\bibitem {BE} \textsc{Bertacchi, D.} (2006). Asymptotic behaviour of the
simple random walk on the 2-dimensional comb. \textit{Electron. J.
Probab.} \textbf{11} 1184--1203.

\bibitem{BZ}
\textsc{Bertacchi, D. and Zucca, F.} (2003). Uniform asymptotic
estimates of transition probabilities on combs. \textit{J. Aust. Math.
Soc.} \textbf{75} 325--353.

\bibitem{CCFR08}
\textsc{Cs\'aki, E., Cs\"org\H o, M., F\"oldes, A. and R\'ev\'esz,
P.} (2009). Strong limit theorems for a simple random walk on the
2-dimensional comb. \textit{Electr. J. Probab.} \textbf{14} 2371--2390.

\bibitem{CCFR11}
\textsc{Cs\'aki, E., Cs\"org\H o, M., F\"oldes, A. and R\'ev\'esz,
P.} (2011). On the local time of random walk on the 2-dimensional comb.
\textit{Stochastic Process. Appl.} \textbf{121} 1290--1314.

\bibitem{CCFR13}
\textsc{Cs\'aki, E., Cs\"org\H o, M., F\"oldes, A. and R\'ev\'esz,
P.} (2013). Strong limit theorems for anisotropic random walks on ${\IZ}^2$.
\textit{Periodica Math. Hungar.} \textbf{67} 71--94.

\bibitem{CFR15} 
\textsc{Cs\'aki, E., F\"oldes, A. and R\'ev\'esz, P.} (2015). Some results 
and problems for anisotropic random walks on the plane. D. Dawson {\sl et al.}
(eds), \textit{Asymptotic Laws and Methods in Stochastics}, Fields Institute 
Communications \textbf{76} 55--75. Springer Science + Business Media New
York.

\bibitem{CM}
\textsc{Cs\"org\H o, M. and Martsynyuk, Y.V.} (2011). Functional central limit
theorems for self-normalized least squares processes in regression with 
possibly infinite variance data. \textit{Stochastic Process. Appl.} 
\textbf{121} 2925--2953.

\bibitem{CR79}
\textsc{Cs\"org\H o, M. and R\'ev\'esz, P.} (1979). How big are the increments
of a Wiener process? \textit{Ann. Probab.} \textbf{7} 731--737. 

\bibitem{CR81}
\textsc{Cs\"org\H o, M. and R\'ev\'esz, P.} (1981). \textit{Strong
Approximation in Probability and Statistics}. Academic, New York.

\bibitem{DH}
\textsc{den Hollander, F.} (1994). On three conjectures by K. Shuler.
\textit{J. Statist. Physics} \textbf{75} 891--918.

\bibitem{DE}
\textsc{Dvoretzky, A. and Erd\H os, P.} (1951). Some problems on random
walk in space. \textit{Proc. Second Berkeley Symposium}, pp. 353--367.

\bibitem{ET}
\textsc{Erd\H os, P. and Taylor, S.J.} (1960). Some problems concerning
the structure of random walk paths. \textit{Acta Math. Acad. Sci.
Hungar.} \textbf{11} 137--162.

\bibitem{HH}
\textsc{Hall, P. and Heyde, C.C.} (1980). \textit{Martingale Limit Theory
and its Application.} Academic, New York.

\bibitem{H}
\textsc{Heyde, C.C.} (1982). On the asymptotic behaviour of random walks
on an anisotropic lattice. \textit{J. Statist. Physics} \textbf{27}
721--730.

\bibitem{H93}
\textsc{Heyde, C.C.} (1993). Asymptotics for two-dimensional anisotropic
random walks. In: \textit{Stochastic Processes.} Springer, New York,
pp. 125--130.

\bibitem{KMT}
\textsc{Koml\'os, J., Major, P. and Tusn\'ady, G.} (1975).
An approximation of partial sums of independent rv's and the sample df.
I. \textit{Z. Wahrsch. verw. Gebiete} \textbf{32} 111--131.

\bibitem{RE}
\textsc{R\'ev\'esz, P.} (2013).
\textit{Random Walk in Random and Non-Random Environments}, 3rd ed.
World Scientific, Singapore.

\bibitem{SH}
\textsc{Shuler, K.E.} (1979). Random walks on sparsely periodic and random
lattices I. \textit{Physica A} \textbf{95} 12--34.

\bibitem{ST}
\textsc{Strassen, V.} (1964).
An invariance principle for the law of the iterated logarithm.
\textit{Z. Wahrsch. verw. Gebiete} \textbf{3} 211--226.

\bibitem{ST67}
\textsc{Strassen, V.} (1967). Almost sure behaviour of sums of independent 
random variables and martingales. \textit{Proceedings of the Fifth Berkeley
Symposium of Mathematical Statistics and Probability,} \textbf{2}, 315--343.
University of California Press, Berkeley. 

\bibitem{WH}
\textsc{Weiss, G.H. and Havlin, S.} (1986). Some properties of a random
walk on a comb structure. \textit{Physica A} \textbf{134} 474--482.

\end{thebibliography}
\end{document}